\documentclass[final,leqno,onefignum,onetabnum]{amsart}

\usepackage{moreverb}
\usepackage{lineno}%,hyperref}
\modulolinenumbers[5]
\usepackage{array,latexsym,pst-node,amsmath,amssymb}

\newcommand{\mystrut}{\vrule height9.5pt depth1.5pt width0pt}
\newcommand{\tab}{\par\noindent\mystrut}
\newcommand{\tabb}{\tab\hskip1.5em}

\newcommand{\alg}[1]{\par\noindent\mystrut\ignorespaces\hbox to\textwidth{#1\hfill}}
\newcommand{\algt}[1]{\par\noindent\mystrut\hbox to\textwidth{\ignorespaces\hskip1.5em#1\hfill}}
\newcommand{\algtt}[1]{\par\noindent\mystrut\hbox to\textwidth{\ignorespaces\hskip3.5em#1\hfill}}
\newcommand{\algttt}[1]{\par\noindent\mystrut\hbox to\textwidth{\ignorespaces\hskip5.5em#1\hfill}}
\newcommand{\algtttt}[1]{\par\noindent\mystrut\hbox to\textwidth{\ignorespaces\hskip7.5em#1\hfill}}
\newcommand{\algttttt}[1]{\par\noindent\mystrut\hbox to\textwidth{\ignorespaces\hskip9.5em#1\hfill}}
\newcommand{\clg}[2]{\par\noindent\mystrut\hbox to\textwidth{\ignorespaces#1\hfill[#2]}}
\newcommand{\clgt}[2]{\par\noindent\mystrut\hbox to\textwidth{\ignorespaces\hskip1.5em#1\hfill[#2]}}
\newcommand{\clgtt}[2]{\par\noindent\mystrut\hbox to\textwidth{\ignorespaces\hskip3.5em#1\hfill[#2]}}
\newcommand{\clgttt}[2]{\par\noindent\mystrut\hbox to\textwidth{\ignorespaces\hskip5.5em#1\hfill[#2]}}
\newcommand{\clgtttt}[2]{\par\noindent\mystrut\hbox to\textwidth{\ignorespaces\hskip7.5em#1\hfill[#2]}}
\newcommand{\clgttttt}[2]{\par\noindent\mystrut\hbox to\textwidth{\ignorespaces\hskip9.5em#1\hfill[#2]}}

\newcommand{\FOR}{\textbf{for}\hskip2pt}

\newcommand{\END}{\textbf{end}}

\usepackage[norelsize]{algorithm2e}  %FIX FOR arxiv

\usepackage{latexsym}
\usepackage{color}
\usepackage{multirow}
\usepackage{tabularx}

\usepackage[dvips,colorlinks,bookmarksopen,bookmarksnumbered,citecolor=red,urlcolor=red]{hyperref}

\newcommand\BibTeX{{\rmfamily B\kern-.05em \textsc{i\kern-.025em b}\kern-.08em
T\kern-.1667em\lower.7ex\hbox{E}\kern-.125emX}}

\newcolumntype{x}[1]{%
&gt;{\raggedleft\arraybackslash}p{#1}}%

\title{Compact Representation of the Full Broyden Class of Quasi-Newton Updates}

\author{Omar DeGuchy}
\email{odeguchy@ucmerced.edu}
\address{School of Natural Sciences, University of California,
    Merced, 5200 N. Lake Road, Merced, CA 95343}
\author{Jennifer B. Erway}
\email{erwayjb@wfu.edu}
\address{Department of Mathematics, PO Box 7388,
    Wake Forest University, Winston-Salem, NC 27109}
\author{Roummel F. Marcia}
\email{rmarcia@ucmerced.edu}
\address{School of Natural Sciences, University of California,
    Merced, 5200 N. Lake Road, Merced, CA 95343}
\thanks{Research supported in part by NSF grants CMMI-1334042 and CMMI-1333326.}

\begin{document}

\newcommand{\subject}{\mathop{\mathrm{subject\ to}}}  % Need \usepackage{amssymb}                               
\newcommand{\minimize}[1]{{\displaystyle\minim_{#1}}}

\newcommand{\minim}{\mathop{\mathrm{minimize}}}
\newcommand{\minimizearg}[2]{{\displaystyle\minim_{#1 \in \Re^{#2}}}}
\newcommand{\maximizearg}[2]{{\displaystyle\maxim_{#1 \in \Re^{#2}}}}
\newcommand{\minimizetwoarg}[4]{{\displaystyle\minim_{#1 \in \Re^{#2},%                                                
                                                  #3 \in \Re^{#4}}}}
\newcommand{\maximize}[1]{{\displaystyle\maxim_{#1}}}
\newcommand{\maxim}{\mathop{\mathrm{maximize}}}
\renewcommand{\maximize}[1]{{\displaystyle\maxim_{#1}}}
\newcommand{\MATLAB}{{\small MATLAB}}
\newcommand{\SMW}{{\small SMW}}
\newcommand{\QR}{{\small QR}}
\newcommand{\BFGS}{{\small BFGS}}
\newcommand{\DFP}{{\small DFP}}
\newcommand{\SR}{{\small SR1}}
\newcommand{\LSR}{{\small L-SR1}}
\newcommand{\LBFGS}{{\small L-BFGS}}
\newcommand{\mgap}{\;\;}
\newcommand{\bgap}{\;\;\;}
\newcommand{\defined}{\mathop{\,{\scriptstyle\stackrel{\triangle}{=}}}\,}
\def\wM{\widetilde{M}}
\def\wPsi{\widetilde{\Psi}}

\makeatletter
\newcommand{\infimum}[1]{{\displaystyle\infim_{#1}}}
\newcommand{\infim}{\mathop{\operator@font{inf}}}
\newcommand{\todo}[1]{{\textcolor{red}{[#1]}}}
\newcommand{\rfm}[1]{{\textcolor{blue}{#1}}}
\newcommand\XOR{\mathbin{\char`\^}}

\newcommand{\algbox}{%
  \vbox{\hrule height0.3pt\hbox{\vrule height2pt%
      width0.3pt\hskip\textwidth\hskip-0.6pt\vrule width0.3pt}\hrule height0.3pt}}

\newcommand{\AlgBegin}{\vspace{\smallskipamount}\noindent\algbox\par\vspace{\smallskipamount}\par\noindent}
\newcommand{\AlgEnd}{\par\noindent\algbox\vspace{\smallskipamount}}

\newcounter{pseudocode}[section]
\def\thepseudocode{\thesection.\arabic{pseudocode}}
\newenvironment{pseudocode}[2]%
        {%
        \refstepcounter{pseudocode}%
          \AlgBegin %
               {{\bfseries Algorithm \thepseudocode.}\rule[-1.25pt]{0pt}{10pt}#1}%
        #2}%
           {\AlgEnd}

\newcounter{Pseudocode}[section]
\def\thePseudocode{\thesection.\arabic{Pseudocode}}
\newenvironment{Pseudocode}[2]%
        {%
        \refstepcounter{Pseudocode}%
          \AlgBegin %
               {{\bfseries #1.}\rule[-1.25pt]{0pt}{10pt}}%
        #2}%
           {\AlgEnd}

\def\tnu{\tilde{\nu}}
\makeatother

\maketitle

\begin{abstract}
In this paper, we present the compact representation for matrices
  belonging to the the Broyden class of quasi-Newton updates, where each
  update may be either rank-one or rank-two.  This work extends previous
  results solely for the \emph{restricted} Broyden class of rank-two
  updates.  In this article, it is not assumed the same Broyden update is
  used each iteration; rather, different members of the Broyden class may
  be used each iteration.  Numerical experiments suggest that a practical
  implementation of the compact representation is able to accurately
  represent matrices belonging to the Broyden class of updates.
  Furthermore, we demonstrate how to compute the compact representation for
  the inverse of these matrices, as well as a practical algorithm for solving
  linear systems with members of the Broyden class of updates.  We
  demonstrate through numerical experiments that the proposed linear solver
  is able to efficiently solve linear systems with members of the
  Broyden class of matrices to high accuracy. As an immediate consequence
  of this work, it is now possible to efficiently compute the eigenvalues
  of any limited-memory member of the Broyden class of matrices, allowing
  for the computation of condition numbers and the ability perform
  sensitivity analysis.
\end{abstract}


\pagestyle{myheadings}
\thispagestyle{plain}
\markboth{DeGuchy \emph{et~al.}}{Compact Representation of the Full Broyden Class of Quasi-Newton Updates}

%%%%%%%%%%%%%%%%%%%%%%%%%%%%%%%%%%%%%%%%%%%%%%%%%%%%%%%%%%
\section{Introduction}
\vspace{-2pt}
Quasi-Newton methods for minimizing a continuously
differentiable function $f:\Re^n\rightarrow \Re$ generate a
sequence of iterates $\{x_k\}$ such that $f$ is strictly decreasing 
at each iterate.  Crucially, at each iteration a quasi-Newton matrix
is used to approximate $\nabla^2f(x_k)$ that is assumed to
be either too computationally expensive to compute or unavailable.
The approximation to the Hessian is updated each iteration using the 
most recently-computed iterate $x_{k+1}$ by defining
a new \emph{quasi-Newton pair} $(s_k, y_k)$ given by
$$
	s_k\defined x_{k+1}-x_k \quad 
	\text{and} \quad 
	y_k\defined \nabla f(x_{k+1})-\nabla f(x_k).
$$
The quasi-Newton Broyden family of updates is given by
\begin{equation}\label{eqn-1param}
% B_{k+1}^\phi =B_k^\phi  -  \frac{1}{s_k ^TB_k^\phi s_k }B_k^\phi s_k s_k ^TB_k^\phi_k +  
%  	\frac{1}{y_k^Ts_k }y_ky_k^T +
%         \phi (s_k^T B_k^\phi s_k )w_kw_k^T,
 B_{k+1} =B_k  -  \frac{1}{s_k ^TB_k s_k }B_k s_k s_k ^TB_k +  
  	\frac{1}{y_k^Ts_k }y_ky_k^T +
         \phi_k (s_k^T B_k s_k )w_kw_k^T,
\end{equation}
 where $\phi_k \in \Re$ and
$$
	w_k = \frac{y_k}{y_k^Ts_k} - \frac{B_k s_k}{s_k^TB_k s_k}.
$$
For $\phi_k \in [0,1]$, $B_{k+1}$ is said to be in the \textsl{restricted}
or {convex} Broyden class of updates.  Setting $\phi=0$ gives the
Broyden-Fletcher-Goldfarb-Shanno (\BFGS) update, arguably the most
widely-used symmetric positive-definite update and a member of the
restricted Broyden class. For $\phi\not\in [0,1]$, the sequence of
quasi-Newton matrices generated by this update is not guaranteed to be
positive definite.  The most well-known update not in the restricted
Broyden class is the symmetric rank-one (\SR1) update, which is obtained by
setting $\phi_k =(s_k^Ty_k)/(s_k^Ty_k - s_k^TB_ks_k)$.

Recently, there has been renewed interest in the entire Broyden class of updates,
and in particular, in negative values of $\phi$.  Research has
shown that negative values of $\phi$ are desirable~\cite{byrd1992behavior}
and under some conditions, quasi-Newton methods based on negative values of
$\phi$ exhibit superlinear convergence
rates~\cite{byrd1992behavior,preconvex}.  There has also been empirical
evidence that $\phi<0$ may lead to more efficient algorithms than
\BFGS~\cite{preconvex,liu2007statistical}.

\medskip

In this paper, we present the compact representation for the full Broyden
class of quasi-Newton matrices, allowing $\phi$ to be negative and to
change each iteration.  We also demonstrate how to efficiently solve linear
systems with any member of the Broyden class using the compact
representation of its inverse.  This paper can be viewed as an extension of
the results found in~\cite{ErwM15,ErwM_CptInv}, which presented the compact
representation for members of the restricted Broyden class and their
inverses, as well as a practical method for solving linear systems
involving only restricted Broyden class matrices (i.e., $\phi\in [0,1]$).

One important application of the compact representation is the ability to
efficiently compute the eigenvalues and a partial eigenbasis when the
number of stored pairs is small~\cite{ErwM15}, which is the case in
large-scale optimization with so-called \emph{limited-memory} quasi-Newton
updates.  In this setting, only the most recently-computed $M$ quasi-Newton
pairs $\{(s_k,y_k)\}$, $k=0,1,\ldots, M-1$, are stored and used to update
$B_{k+1}$ using the recursive application of (\ref{eqn-1param}).
Typically, in large-scale applications $M\le 10$ regardless of $n$, i.e.,
$M\ll n$ (see, e.g., ~\cite{ByrNS94}).  With the eigenvalues it is now
possible to compute condition numbers, compute singular values, and perform
sensitivity analysis.

This paper is organized in seven sections.  In the second section, we
review the compact formulation for the restricted Broyden class of
updates ($\phi\in [0,1]$) as well as overview the efficient computation of their eigenvalues.
The main result of the paper is in Section 3 where the compact
representation is given for the entire Broyden class of updates that allows
for $\phi$ to change each update.  In this section, we also present a
practical iterative method to compute the compact representation.  In
Section 4, we show how to perform linear solves with any member of the
Broyden class using the compact representation of their inverse.
Numerical experiments are reported in Section 5. Finally, Section 6
contains concluding remarks, and Section 7 includes acknowledgements for
this work.

\subsection{Notation and assumptions}
Throughout this paper, we make use of the following matrices:
\begin{eqnarray}\label{eqn-SY}
	S_k &\defined& \big ( \ s_0 \ \ s_1 \ \ s_2 \ \ \cdots \ \ s_{k} \ \big ) \ \in \ \Re^{n \times (k+1)}, \\
	Y_k &\defined& \big ( \ y_0 \ \ y_1 \ \ y_2 \ \ \cdots \ \ y_{k} \ \big ) \ \in \ \Re^{n \times (k+1)}.
\end{eqnarray}
Furthermore, we make use of the following decomposition of $S_k^TY_k \in\Re^{(k+1) \times (k+1)}$:
\begin{equation}\label{eqn-STY}
	S_k^TY_k =   L_k + D_k + R_k,
\end{equation}
where $L_k$ is strictly lower triangular, $D_k$ is diagonal, and $R_k$ is
strictly upper triangular.  
%We assume that $s_i^Ty_i > 0$ for all $i = 0,
%1, \dots, k$.\od{Are we still assuming the previous statement?} 
We assume that the matrix $B_k$ is
nonsingular for each $k$.
Finally, throughout the manuscript, $I_j$ denotes the $j\times j$ identity matrix.

\vspace{-6pt}
\section{Compact representation for the restricted Broyden class} 
\vspace{-2pt}

Compact representations of matrices from the Broyden class of updates were
first described by Byrd et al~\cite{ByrNS94} as matrix decompositions of
the form
$$B_{k+1}=B_0+\Psi_k M_k\Psi_k^T,$$ where 
$\Psi_k\in\Re^{n\times l}$, $M_k\in\Re^{l\times l}$, and $B_0$ is the initial
matrix.  The size of $l$ depends on the rank of the update; in the case
of a rank-two update, $l=2(k+1)$, and in the case of a rank-one update,
$l=k+1$.  In the case of the \BFGS{} update (i.e., $\phi=0$), $\Psi_k$ and $M_k$ are given in~\cite{ByrNS94}:
\begin{equation}\label{eqn-alt-form2}
	\Psi_k\defined \begin{pmatrix} B_0S_k & Y_k \end{pmatrix} \quad \text{and} \quad
	M_k \defined
	-\begin{pmatrix}
	S_k^TB_0S_k & L_k \\
	L_k^T & -D_k
	\end{pmatrix}^{-1},
\end{equation}
where $S_k$ and $Y_k$ are defined in (\ref{eqn-SY}).
In~\cite{ErwM15}, we presented the compact representation for any matrix in
the restricted Broyden class (i.e, $\phi\in [0,1]$); in particular, for any
matrix in the restricted Broyden class,
$$
	\Psi_k\defined\begin{pmatrix} B_0S_k & Y_k \end{pmatrix} \quad \text{and} \quad
M_k=	\begin{pmatrix}
		-S_k^TB_0S_k + \phi \Lambda_k  & -L_k + \phi \Lambda_k \\
		-L_k^T + \phi \Lambda_k & \ \ D_k + \phi \Lambda_k
	\end{pmatrix}^{-1},$$
where $L_k$ and  $D_k$ are given in (\ref{eqn-STY}) and
$\Lambda_k \in \Re^{(k+1) \times (k+1)}$ is
the diagonal matrix $\Lambda_k=\text{diag}(\lambda_i),$ ($0\le i \le k$),
given by
\begin{equation}\label{eq:Lambda}
			\lambda_{i}\defined
		%	\frac{1}{\displaystyle 
		\left(	-\frac{1-\phi}{s_i^TB_is_i}
			-\frac{\phi}{s_i^Ty_i}\right)^{-1}.%}.
%\quad	\ \text{for $0 \le i \le k$}.
\end{equation}

To our knowledge, the only compact formulation known for a member of the Broyden
class of updates outside the restricted class is for
an \SR1{} matrix.
As with the \BFGS{} case, it is also given in~\cite{ByrNS94};
in particular,
\begin{equation*}
	\Psi_k \ = \ 
       	Y_k  - B_0S_k \quad \text{and} \quad
        M_k \ = \ (D_k + L_k + L_k^T - S_k^TB_0S_k)^{-1}.
\end{equation*}
Notice that $\Psi_k$ in the compact representation for \SR1{} matrices is
half the size of that of $\Psi_k$ for the rank-two updates.

\subsection{Applications of the compact representation}
In this section, we briefly review how the eigenvalues of any quasi-Newton
matrix that exhibits a compact representation can be efficiently computed.
The first method to compute eigenvalues of limited-memory quasi-Newton
matrices was proposed by Lu~\cite{Lu96}.  This method makes use of the
singular value decomposition and an eigendecomposition of small matrices.
An alternative approach, first described by Burdakov et
al.~\cite{Burdakov13}, uses the \QR{} factorization in lieu of the
singular value decomposition.  An overview of the method found
in~\cite{Burdakov13} follows below.  For this section, we assume $k$ is
small, as in the case of limited-memory quasi-Newton matrices; moreover, we
assume $\Psi_k\in\Re^{n\times l}$ is full rank, where $l$ is either
$l=2(k+1)$ or $l=k+1$.  Finally, we assume $B_0=\gamma I$, where
$\gamma\in\Re$.

Let $QR$ be the ``thin'' \QR{} decomposition of $\Psi_k$, where $Q\in \Re^{n
  \times l}$ has orthonormal columns and $R\in \Re^{l \times l}$ is upper
triangular (see, e.g., \cite{GVL96}).
Then, \begin{equation*}
	B_{k+1}  =  B_0  + \Psi_k M_k \Psi_k^T =
B_0 + Q R M_kR^TQ^T.
\end{equation*}
The matrix
$RM_kR^T$ is a real symmetric $l\times l$ matrix, whose
spectral decomposition can be explicitly computed since $l$ is small.  Letting 
$VDV^T$ be its spectral decomposition gives that
\begin{equation} \label{eqn-eigs}
	B_{k+1}
	= B_0 + Q V D V^T Q^T 
	= \gamma I + Q V D V^T Q^T 
	= Q V ( \gamma I+ \hat{D} ) V^T Q^T,
\end{equation}
where $\hat{D}$ is a diagonal matrix whose leading $l\times l$ block is $D$ while the rest of
the matrix is zeros.  Thus, the spectral decomposition of $B_{k+1}$ is
given by (\ref{eqn-eigs}).  (Note that in practice, the matrices $Q$ and
$V$ in (\ref{eqn-eigs}) are not stored.)  Note that the matrix $B_{k+1}$ has an
eigenvalue of $\gamma$ with multiplicity $n-l$ and $l$ eigenvalues given by
$\gamma + d_i$, where $1 \le i \le l$ and $D_{i,i}\defined d_i$.  It also turns out that it is also
possible to efficiently compute the eigenvectors associated with the
nontrivial eigenvalues and only one eigenvector associated with the trivial
eigenvalue $\gamma$.  (For more details, see~\cite{Burdakov13, ErwM15}.)

\medskip

Generally speaking, computing the eigenvalues of $B_{k+1}$ directly is an
$O(n^3)$ process.  In contrast, the above decomposition requires the \QR{}
factorization of $\Psi_k$ and the eigendecomposition of $RM_kR^T$,
requiring $O(nl^2)$ flops and $O(l^3)$ flops, respectively.  Since $l\ll
n$, the proposed method's runtimes should increase only linearly with $n$.
(For some details regarding updating the (full) \QR{} factorization after a
new quasi-Newton pair is computed, see~\cite{ErwM15}.)  This efficient
computation of eigenvalues and a partial eigenbasis appears in new methods
for large-scale
optimization~\cite{Burdakov13,BruEM16,Lasith_isita,Shelby,Brust_shape}.

The compact representation is also useful for solving linear systems with
quasi-Newton matrices.  In~\cite{BurWX96}, Burke et al. use the compact
formulation of a \BFGS{} matrix to solve a linear system involving a
diagonally-shifted \BFGS{} matrix. In~\cite{ErwM_CptInv}, the compact
representation for the inverse of any member in the restricted Broyden class
is given as well as a practical method to solve linear systems involving
these matrices using this representation.

\vspace{-6pt}
\section{Compact Representation for any member of the Broyden class}
\vspace{-2pt}

The main result for this section is a theorem giving the compact
representation for any member of the Broyden class.  The representation
allows $\phi$ to change each iteration and to be negative.  In this
section, we also present a practical algorithm for computing the compact
representation.

\medskip

We begin by observing that $B_{k+1}$ in \eqref{eqn-1param} can be written as
 \begin{equation}\label{eqn-compact}
	B_{k+1}
	=
	B_k 
	+ 
	\left ( B_k s_k \ \ y_k \right )
	O_k
	\begin{pmatrix}
		(B_ks_k)^T \\
		y_k^T
	\end{pmatrix},
\end{equation}
where
\begin{equation}\label{eq:Ok}
	O_k 
	=
	\begin{pmatrix}
		\displaystyle - \frac{(1-\phi_k)}{s_k^TB_ks_k} 
		&  \displaystyle - \frac{\phi_k}{y_k^Ts_k} \\
		\displaystyle -\frac{\phi_k}{y_k^Ts_k} 
		& \displaystyle  \left ( 1 + \phi_k \frac{s_k^TB_ks_k}{y_k^Ts_k}\right ) \frac{1}{y_k^Ts_k} 
	\end{pmatrix}.
\end{equation}

% We begin by considering the \SR1{} update, which is the unique rank-one
% update in the Broyden class.  This update is obtained by picking 
% $\phi$ to be $\phi_k^{\text{\SR1}}
% \defined (s_k^Ty_k)/(s_k^Ty_k - s_k^TB_ks_k)$.  Note that $\phi^{\text{SR1}}_k\ne 0$ since
% $s_k^Ty_k \ne 0$; moreover, since $B_k$ is assumed to be invertible,
% $s_k^TB_ks_k \ne 0$, and thus, $\phi_k^{\SR1} \ne 1$.  Before stating the 
% main theorem, two lemmmas are given.
\noindent
We now state two lemmas about $O_k$; specifically, we provide the condition
for $O_k$ when is singular as well as its inverse when it is
nonsingular. \bigskip

\noindent 
\textbf{Lemma 1.} \textsl{The $2 \times 2$ matrix $O_k$ is singular if and only if $\phi_k= (s_k^Ty_k)/(s_k^Ty_k - s_k^TB_ks_k)$.}

\bigskip

\noindent \textbf{Proof.}
%We first consider two special cases: $\phi_k=0$ and $
%\phi_k=1$.  Notice that if $\phi_k = 0$ then $\phi_k \ne 
%(s_k^Ty_k)/(s_k^Ty_k - s_k^TB_ks_k)$ since 
%$y_k^Ts_k \ne 0$.  Also, if $\phi_k = 1$, then $\phi_k \ne
%(s_k^Ty_k)/(s_k^Ty_k - s_k^TB_ks_k)$ since $B_k$ is assumed to be nonsingular,
%and thus,
%$s_k^TB_ks_k \ne 0$. Using (\ref{eq:Ok}) it can be seen that if $\phi_k \in \{0,1\}$, then $O_k$ is nonsingular.
%Thus, we have shown that the lemma holds for $\phi_k = 0$ and $\phi_k=1$.  
%We now consider $\phi_k\not\in\{0,1\}.$
% Let
% \begin{eqnarray*}
%         \alpha_k = -\frac{1-\phi_k}{s_k^TB_ks_k}, \quad
%         \beta_k  = - \frac{\phi_k}{y_k^Ts_k}, \quad \text{and} \quad 
%         \delta_k  =   \left ( 1 + \phi_k \frac{s_k^TB_ks_k}{y_k^Ts_k}\right ) \frac{1}{y_k^Ts_k}.
% \end{eqnarray*}
% Then the determinant, $\eta_k$, of $O_k$  can be written as 
% \begin{eqnarray}\label{eqn-etak}
% 	\eta_k
% 	=
% 	\alpha_k\delta_k - \beta_k^2
% 	=	
% 	\frac{1}{y_k^Ts_k}(\alpha_k + \beta_k).
% \end{eqnarray}
% Thus, $\eta_k = 0$, i.e., $O_k$ is singular,  if and only if $\alpha_k + \beta_k = 0$, i.e.,
The determinant of $O_k$ is given by
%\begin{eqnarray*}\det(O_k) &=&
%- \frac{(1-\phi_k)}{s_k^TB_ks_k} \left ( 1 + \phi_k \frac{s_k^TB_ks_k}{y_k^Ts_k}\right ) \frac{1}{y_k^Ts_k} -  \frac{\phi^2_k}{\left(y_k^Ts_k\right)^2}\\
%& = &  \frac{1}{y_k^Ts_k} \left( \frac{(1-\phi_k)}{s_k^TB_ks_k}- \frac{\phi_k}{y_k^Ts_k}\right).
%\end{eqnarray*}
\begin{eqnarray*}\det(O_k) &=&
	- \frac{(1-\phi_k)}{s_k^TB_ks_k} \left ( 1 + \phi_k \frac{s_k^TB_ks_k}{y_k^Ts_k}\right ) \frac{1}{y_k^Ts_k} -  \frac{\phi^2_k}{\left(y_k^Ts_k\right)^2}\\
	& = &  \frac{1}{y_k^Ts_k} \left( -\frac{(1-\phi_k)}{s_k^TB_ks_k}- \frac{\phi_k}{y_k^Ts_k}\right).
\end{eqnarray*}
Thus, $O_k$ is singular if and only if 
%\begin{equation}\label{eqn-singO}
%	- \frac{(1-\phi_k)}{s_k^TB_ks_k} + \frac{\phi_k}{y_k^Ts_k} = 0;
%\end{equation}
\begin{equation}\label{eqn-singO}
- \frac{(1-\phi_k)}{s_k^TB_ks_k} - \frac{\phi_k}{y_k^Ts_k} = 0; 
\end{equation}
in other words, $\phi_k = (y_k^Ts_k)/(y_k^Ts_k - s_k^TB_ks_k)$.
$\square$
	
Lemma 1 states that $O_k$ is singular if and only if the \SR1{} update is
used.  Special care will given to the \SR1{} case, since unlike other
members of the Broyden class, this is a rank-one update.  For the duration
of this manuscript, we let $\phi_k^{\SR1}\defined
(y_k^Ts_k)/(y_k^Ts_k - s_k^TB_ks_k)$.  For $\phi_k\ne \phi_k^{\SR1}$,
$O_k$ is invertible and its inverse is given in Lemma 2.  This result can
be derived by using the formula for the inverse of a $2\times 2$ matrix.

\bigskip

\noindent \textbf{Lemma 2.} \textsl{If $O_k$ is invertible, then}
$$
	O_k^{-1}
	=
	\begin{pmatrix}
		\displaystyle -s_k^TB_ks_k +  \frac{\phi_k}{\alpha_k + \beta_k} 
		& 
		\displaystyle \frac{\phi_k}{\alpha_k + \beta_k} 
		\\
		\displaystyle \frac{\phi_k}{\alpha_k + \beta_k} 
		& 
		\displaystyle y_k^Ts_k +  \frac{\phi_k}{\alpha_k + \beta_k}
	\end{pmatrix},
$$
where $\alpha_k = -(1-\phi_k)/(s_k^TB_ks_k)$ and
$\beta_k  = - \phi_k/(y_k^Ts_k).$

\bigskip

\bigskip

We now state the main theorem of this manuscript that presents the compact
representation for any member of the Broyden class, while allowing the
parameter $\phi$ to vary at each iteration.  After proving this theorem, we
discuss several aspects of this compact representation as well as the key
differences between the compact representation for the Broyden class
of matrices (Theorem 1) and the compact representation of the \emph{restricted}
Broyden class reviewed in Section 2.

\bigskip

\noindent 
\textbf{Theorem 1.} 
\textsl{
Let $\Psi_k = (B_0S_k \ \ Y_k ) \in \Re^{n \times 2(k+1)}$.
Let $\Pi_k \in \Re^{2(k+1) \times 2(k+1)}$ be the permutation matrix
\begin{equation}\label{eq:Pi_k}
	{\Pi_k}
	=
	\begin{pmatrix}
		I_k & 0 & 0 & 0 \\
		0 & 0 & I_k & 0 \\
		0 & 1 & 0 & 0 \\
		0 & 0 & 0 & 1
	\end{pmatrix},  
\end{equation}
with $\Pi_0\defined I_2$.
Additionally, let $\Xi_k$ be defined recursively as
\begin{equation}\label{eq:Xi_k}
	{\Xi_k}
	=
		\begin{pmatrix}
			\Pi_{k-1}^T\Xi_{k-1} & 0 \\
			0 & E_{k}
		\end{pmatrix},
		\quad \text{where } \quad
		E_k = 
		\begin{cases} \ \ 
                  (-1 \ \ 1)^T 
		& \text{if  $\phi_k = \phi_k^{\SR1}$}\\
%		\\[.5cm]
		I_2 & \text{otherwise,}
		\end{cases}
%		\begin{cases} \ \ 
%		\begin{pmatrix}
%			-1 \\
%			\ \ \  \! 1 
%		\end{pmatrix}
%		& \text{if  $\phi_k = \phi_k^{\SR1}$}
%		\\[.5cm]
%		\begin{pmatrix}
%			1 & 0 \\
%			0 & 1
%		\end{pmatrix}
%		& \text{otherwise}
%		\end{cases},
\end{equation}
with $\Xi_0 \defined E_0$.
Finally, let $\Gamma_k \in \Re^{(k+1) \times (k+1)}$ be
a diagonal matrix such that
\begin{equation}\label{eq:Gamma}
	\Gamma_k =  \underset{0 \le j \le k}{\text{diag}} \big ( \gamma_j \big ),
	\qquad \text{where \ } 
			\gamma_j = 
			\begin{cases}
\phi_j\left(	-\frac{1-\phi_j}{s_j^TB_js_j}
			-\frac{\phi_j}{s_j^Ty_j}\right)^{-1}
		%	\frac{\displaystyle \phi_j}{\displaystyle 
			%-\frac{1-\phi_j}{s_j^TB_js_j}
			%-\frac{\phi_j}{s_j^Ty_j}}
			& \text{if} \  \phi_j \ne \phi_j^{\SR1}
			\\
			\qquad \quad \ 0 
			& \text{otherwise.}
			\end{cases}
\end{equation}
If $B_{k+1}$ is a member of the Broyden class of updates, then
\begin{equation}\label{eq:compact}
	B_{k+1} = B_0 + \widehat{\Psi}_k  \widehat{M}_k \widehat{\Psi}_k^T,
\end{equation} 
where
\begin{equation}\label{eqn-G} 
	\widehat{M}_k=	
	\left (
	\Xi_k^T \Pi_k
	\begin{pmatrix}
		-S_k^TB_0S_k + \Gamma_k  & -L_k + \Gamma_k \\
		-L_k^T + \Gamma_k & \ \ D_k + \Gamma_k
	\end{pmatrix}
	\Pi_k^T
	\Xi_k
	\right )^{-1},
\end{equation}
$L_k$ and $D_k$ are defined in (\ref{eqn-STY}), and
\begin{equation}
\label{eq:Phihat}
	\widehat{\Psi}_k = \Psi_k \Pi_k^T \Xi_k.
\end{equation}
}

\bigskip

\noindent \textbf{Proof.}  This proof is by induction on $k$.
For the base case ($k=0$),  $D_0 = y_0^Ts_0$ with
$L_0 = R_0 = 0$, and $\Gamma_0$ is the scalar $\gamma_0$.  
Thus, $\widehat{M}_0$ defined in (\ref{eqn-G}) reduces to
\begin{equation}\label{eqn-widehatM1}
	\widehat{M}_0 =
	\left (
		\Xi_0^T
		\Pi_0
		\begin{pmatrix}
			-s_0^TB_0s_0 + \gamma_0 & \gamma_0 \\
			\gamma_0 & d_0 + \gamma_0
		\end{pmatrix}
		\Pi_0^T
		\Xi_0
	\right )^{-1}.
\end{equation}
By (\ref{eqn-compact}), $B_1$ is given by
$B_1=B_0+\Psi_0 M_0 \Psi_0^T$ where $\Psi_0 = (B_0s_0 \ \ y_0)$
and  \begin{equation}\label{eq:G0}
	M_0 \defined
		\begin{pmatrix}
			\displaystyle - \frac{(1-\phi_0)}{s_0^TB_0s_0} 
			&  \displaystyle - \frac{\phi_0}{y_0^Ts_0} \\
			\displaystyle -\frac{\phi_0}{y_0^Ts_0} 
			& \displaystyle  \left ( 1 + \phi_0 \frac{s_0^TB_0s_0}{y_0^Ts_0}\right ) \frac{1}{y_0^Ts_0} 
		\end{pmatrix}.
\end{equation}
It remains to show that $M_0 = \Pi_0^T\Xi_0\widehat{M}_0\Xi_0^T\Pi_0$.
Since the initial permutation matrix is defined as $\Pi_0 = I_2$, we only need to show
$ M_0 = \Xi_0 \widehat{M}_0 \Xi_0^T.$
For simplicity, $M_0$ can be written as
\begin{equation}	\label{eqn-widehatM0}
	M_0	=
		\begin{pmatrix}
			\alpha_0 & \beta_0 \\
			\beta_0 & \delta_0
		\end{pmatrix},
\end{equation}
where
\begin{equation}\label{eqn-alphabetadelta}
	\alpha_0 = - \frac{(1-\phi_0)}{s_0^TB_0s_0}, \quad
	\beta_0 = - \frac{\phi_0}{y_0^Ts_0}, \quad \text{and} \quad 
	\delta_0  = \left ( 1 + \phi_0 \frac{s_0^TB_0s_0}{y_0^Ts_0}\right ) \frac{1}{y_0^Ts_0}.
\end{equation}
From Lemma 1, $M_0$ is nonsingular if and only if $\phi_0 \ne
\phi_0^{\SR1}$.  Thus, we consider the following two cases
separately: (a) $\phi_0 = \phi_0^{\SR1}$ and (b) $\phi_0 \ne
\phi_0^{\SR1}$.
\medskip

\noindent \textbf{Case (a):} If $\phi_0 = \phi_0^{\SR1}$,
then $\Xi_0 = E_0 = (-1 \ \ 1)^T$ by \eqref{eq:Xi_k}.
By (\ref{eqn-singO}),
$\alpha_0 + \beta_0 = 0$, and thus,
$M_0$ can be simplified as
$$
	M_0 = 
	\begin{pmatrix}
		-\beta_0 & \ \ \beta_0 \\
		\ \ \beta_0 & -\beta_0
	\end{pmatrix}
	=
	(-\beta_0)
	\begin{pmatrix}
		\ \ 1 & -1\\
		- 1 & \ \ 1
	\end{pmatrix}
	=
	-\beta_0
	\Xi_0 \Xi_0^T.
$$
Finally, since $\phi_0 =  (y_0^Ts_0) / (s_0^Ty_0 - s_0^TB_0s_0)$ and $\beta_0 = -\phi_0/y_0^Ts_0$, then
$$
\widehat{M}_0=
\left (
	\Xi_0^T
	\begin{pmatrix}
		-s_0^TB_0s_0 + \gamma_0 & \gamma_0 \\
		\gamma_0 & s_0^Ty_0 + \gamma_0
	\end{pmatrix}
	\Xi_0
	\right )^{-1}
=	\frac{1}{s_0^Ty_0 - s_0^TB_0s_0} = 
	-\beta_0, 
$$
and thus, $M_0 = \Xi_0 \widehat{M}_0\Xi_0^T$, as desired.

\medskip

\noindent \textbf{Case (b):} If $\phi_0 \ne \phi_0^{\SR1}$, then $M_0$ is nonsingular and $\Xi_0=I_2$.
Thus, it remains to show $M_0 = \widehat{M}_0$.
By Lemma 1, $\alpha_0 + \beta_0 \ne 0$, making $\gamma_0 = \phi_0 / (\alpha_0 + \beta_0)$
 well defined.  By Lemma 2, the inverse of $M_0$ is given by
\begin{equation}\label{eqn-base}
	M_0^{-1}  =
	\begin{pmatrix}
		-s_0^TB_0s_0 + \gamma_0 &  \gamma_0 \\
		\gamma_0 & s_0^Ty_0 + \gamma_0
	\end{pmatrix}
	=
	\widehat{M}_0 ^{-1}.
\end{equation}
Note that the last equality in (\ref{eqn-base}) follows since $\Pi_0=I_2$.
%Thus, the base case is also true when $\phi_k \ne \phi_k^{\SR1}$.
% Note that
% \begin{equation}\label{eq:widehatM0}
% 	\widehat{M}_0 =
% 	\begin{cases}
% 	 -\beta_0
% 	& \text{if $\phi_0 = \phi_0^{\SR1}$}
% 	\\
% 	\begin{pmatrix}
% 		 \alpha_0 & \beta_0 \\
% 		\beta_0 & \delta_0
% 	\end{pmatrix}
% 	& \text{otherwise}
% 	\end{cases}.
% \end{equation}
 
\bigskip
For the induction step, assume 
\begin{equation}\label{eqn-inductB}
	B_{m}=B_0+\widehat{\Psi}_{m-1}  \widehat{M}_{m-1} \widehat{\Psi}_{m-1}^T,
\end{equation} 
where 
$\widehat{M}_{m-1}=	
	\left (
	\Xi_{m-1}^T 
	\Pi_{m-1}
	\Omega_{m-1}
	\Pi_{m-1}^T
	\Xi_{m-1}
	\right )^{-1}$
and 
\begin{equation}\label{eqn-Omega}
	\Omega_{m-1} = 
	\begin{pmatrix}
		-S_{m-1}^TB_0S_{m-1} + \Gamma_{m-1}  & -L_{m-1} + \Gamma_{m-1} \\
		-L_{m-1}^T + \Gamma_{m-1} & \ \ D_{m-1} + \Gamma_{m-1}
	\end{pmatrix}.
\end{equation}
From (\ref{eqn-compact}), we have
\begin{eqnarray}\label{eqn-bnp1}
        B_{m+1} = B_0 \!+\! \widehat{\Psi}_{m-1} \widehat{M}_{m-1} \widehat{\Psi}_{m-1}^T 
        \! + \!
                        \left ( B_ms_m \ \ \ y_m \right )\!
	                \begin{pmatrix}
                                \alpha_m & \beta_m \\
	                        \beta_m & \delta_m
                        \end{pmatrix}\!
             	        \begin{pmatrix}
                                (B_ms_m)^T \\
   	                        y_m^T
                        \end{pmatrix}\!,
\end{eqnarray}
where
\begin{eqnarray*}
        \alpha_m = -\frac{1-\phi_m}{s_m^TB_ms_m}, \quad
        \beta_m  = - \frac{\phi_m}{y_m^Ts_m}, \quad \text{and} \quad 
        \delta_m  =   \left ( 1 + \phi_m \frac{s_m^TB_ms_m}{y_m^Ts_m}\right ) \frac{1}{y_m^Ts_m}.
%      =   \displaystyle - \left ( 1 + (1-\phi_m)\frac{\beta_m}{\alpha_m} \right ) \
%	\frac{\beta_m}{\phi_m}.
\end{eqnarray*}
Multiplying (\ref{eqn-inductB}) by $s_m$ on the right, we obtain
\begin{equation}\label{eqn-block22}
B_ms_m = B_0s_m + \widehat{\Psi}_{m-1}  \widehat{M}_{m-1} \widehat{\Psi}_{m-1}^Ts_m.
\end{equation}
Then, substituting this into (\ref{eqn-bnp1}) yields 
\begin{eqnarray}\label{eqn-bnp12}
	B_{m+1} &=& 
	B_0 \ + \ \widehat{\Psi}_{m-1} \widehat{M}_{m-1} \widehat{\Psi}_{m-1}^T \ + \  \\[.2cm]
	&&
			\left ( B_0s_m\! +\! \widehat{\Psi}_{m-1} p_m \ \ \ y_m \right )\!\!
			\begin{pmatrix}
				\alpha_m & \beta_m \\
				\beta_m & \delta_m
			\end{pmatrix}\!\!
			\begin{pmatrix}
				(B_0s_m \!+\! \widehat{\Psi}_{m-1}p_m)^T \\
				y_m^T
			\end{pmatrix}\!, \nonumber
\end{eqnarray}
where $p_m\defined \widehat{M}_{m-1} \widehat{\Psi}_{m-1}^Ts_m$.
Equivalently,
\begin{eqnarray} 
	B_{m+1}
	&=& B_0 +  \left ( \widehat{\Psi}_{m-1} \ \ B_0s_m \ \ y_m \right )
		\mathcal{M}_{m}
		\begin{pmatrix}
			\widehat{\Psi}_{m-1}^T \\
			(B_0s_m)^T \\
			y_m^T
		\end{pmatrix},
	\label{eqn-3x3}
\end{eqnarray}
where
\begin{equation}\label{eq:Mbar}
\mathcal{M}_m = 
		\begin{pmatrix}
			\widehat{M}_{m-1}+ \alpha_m p_mp_m^T & \alpha_m p_m & \beta_m p_m \\
			\alpha_m p_m^T & \alpha_m & \beta_m \\
			\beta_m  p_m^T & \beta_m & \delta_m
		\end{pmatrix}.
\end{equation}
Note that $\mathcal{M}_m$  has the following decomposition:
\begin{equation}\label{eq:GkP}
%		\begin{pmatrix}
%			M_{m-1}+ \alpha_m p_mp_m^T & \alpha_m p_m & \beta_m p_m\\
%			\alpha_m  p_m ^T & \alpha_m  & \beta_m \\
%			\beta_m  p_m ^T & \beta_m  & \delta_m 
%		\end{pmatrix}
		\mathcal{M}_m
		=
		\begin{pmatrix}
			I & p_m  & 0 \\
			0 & 1 & 0 \\
			0 & 0 & 1
		\end{pmatrix}
		\begin{pmatrix}
			\widehat{M}_{m-1} & 0 & 0 \\
			0 & \alpha_m  & \beta_m  \\
			0 & \beta_m  & \delta_m 
		\end{pmatrix}
		\begin{pmatrix}
			I & 0 & 0 \\
			p_m ^T & 1 & 0 \\
			0 & 0 & 1
		\end{pmatrix}.
\end{equation}

Thus, $\mathcal{M}_m$ is nonsingular if and only if $\alpha_m \delta_m - \beta_m^2 \ne 0$; that is, $\mathcal{M}_m$ is nonsingular if and only if $\phi_m \ne \phi_m^{\SR1}$ (see Lemma 1).
To complete the induction step, we will show that the last term in \eqref{eqn-3x3} is equal to $\widehat{\Psi}_m \widehat{M}_m \widehat{\Psi}_m^T$ by
considering the following two cases separately: (i) $\phi_m = \phi_m^{\SR1}$ and (ii) $\phi_m \ne \phi_m^{\SR1}$.

\bigskip

\noindent \textbf{Case (i):}  If $\phi_m = \phi_m^{\SR1}$,
then by Lemma 1, $\alpha_m = -\beta_m \ne 0$.  Then
\begin{equation}\label{eq:GkPnew}
	\mathcal{M}_m
		=
		\begin{pmatrix}
			I & \ \ 0 \\
			0 & -1 \\
			0 & \ \ 1
		\end{pmatrix}
		\widetilde{\mathcal{M}}_m
		\begin{pmatrix}
			I & \ \ 0 & 0 \\
			0 & -1 & 1
		\end{pmatrix},	
\end{equation}
where
\begin{equation}\label{eqn-Mtildem-i}
	\widetilde{\mathcal{M}}_{m} = 
		\begin{pmatrix}
			\widehat{M}_{m-1} - \beta_mp_mp_m^T &  \beta_mp_m \\
			\beta_m p_m^T & -\beta_m
		\end{pmatrix}.
\end{equation}
We now show that $\widetilde{\mathcal{M}}_m = \widehat{M}_m$.  
By the inductive hypothesis, $\widehat{M}_{m-1}$ is nonsingular. Together
with the fact that $\beta_m\ne 0$, it can be checked directly that
\begin{equation}\label{eqn-minv}
	\widetilde{\mathcal{M}}_m^{-1}
		=
 		\begin{pmatrix}
 			\widehat{M}_{m-1}^{-1} & \widehat{M}_{m-1}^{-1}p_m \\
 			p_m^T\widehat{M}_{m-1}^{-1}& -\beta_m^{-1} + p_m^T\widehat{M}_{m-1}^{-1}p_m
 		\end{pmatrix}.\\
\end{equation}
The (2,2)-entry of $\widetilde{\mathcal{M}}_m^{-1}$ can be simplified
by substituting in for $p_m$ and using the inductive step \eqref{eqn-inductB}:
\begin{eqnarray*}
	-\beta_m^{-1} + p_m^T\widehat{M}_{m-1}^{-1}p_m 
	&=&  
	-\beta_m^{-1} + s_m^T(\widehat{\Psi}_{m-1}\widehat{M}_{m-1}\widehat{\Psi}_{m-1}^T)s_m \\
	&=&
	-\beta_m^{-1} - s_m^TB_0s_m + s_m^TB_m s_m \\
	&=&
	-\beta_m^{-1}  - s_m^TB_0s_m + (1 - \phi_m)\beta_m^{-1} \\
%	&=& 
%	-s_m^TB_0s_m - \phi_m\beta_m^{-1}\\
	&=& -s_m^TB_0s_m + s_m^Ty_m.
\end{eqnarray*}
Substituting this into (\ref{eqn-minv}) and using the inductive hypothesis
gives:
%\begin{eqnarray}
%	\widetilde{\mathcal{M}}_m^{-1} &=&	\begin{pmatrix}
%			\widehat{M}_{m-1}^{-1} & \widehat{\Psi}_{m-1}^Ts_m \\
%			s_m^T\widehat{\Psi}_{m-1} & -s_m^TB_0s_m + s_m^Ty_m
%	\end{pmatrix}	 \nonumber	\\
%	&=&
%	\begin{pmatrix}
%			\left (
%			\Xi_{m-1}^T \Pi_{m-1}
%			\Omega_{m-1} \Pi_{m-1}^T
%			\Xi_{m-1}
%			\right )
%	 		& \Xi_{m-1}^T\Pi_{m-1}\Psi_{m-1}^Ts_m \\
%			s_m^T\Psi_{m-1}\Pi_{m-1}^T\Xi_{m-1} & -s_m^TB_0s_m + s_m^Ty_m
%	\end{pmatrix} \nonumber	 \\	 
%	&=&
%	\begin{pmatrix}
%		\Xi_{m-1}^T\Pi_{m-1}^T & 0 \\
%		0 & 1 
%	\end{pmatrix}
%	\begin{pmatrix}
%	\Omega_{m-1}
%	 & \Psi_{m-1}^Ts_m \\
%	s_m^T\Psi_{m-1}^T & -s_m^TB_0s_m + s_m^Ty_m
%	\end{pmatrix}
%	\begin{pmatrix}
%		\Pi_{m-1}^T\Xi_{m-1} & 0 \\
%		0 & 1 
%	\end{pmatrix}\nonumber\\ 
% 	&=&
%%	\begin{pmatrix}
%%		\Xi_{m-1}^T\Pi_{m-1}^T & \ \ 0 & 0 \\
%%		0 & -1 & 1
%%	\end{pmatrix}
%	\Xi_m^T
%	\begin{pmatrix}
%	\Omega_{m-1}
%	 & -\Psi_{m-1}^Ts_m & 0 \\
%	-s_m^T\Psi_{m-1}^T & -s_m^TB_0s_m + \gamma_m & \gamma_m \\
%	0 & \gamma_m &  s_m^Ty_m + \gamma_m
%	\end{pmatrix}
%	\Xi_m,  \label{eqn-middle}
%\end{eqnarray}
\begin{eqnarray}
\widetilde{\mathcal{M}}_m^{-1} &=&	\begin{pmatrix}
\widehat{M}_{m-1}^{-1} & \widehat{\Psi}_{m-1}^Ts_m \\
s_m^T\widehat{\Psi}_{m-1} & -s_m^TB_0s_m + s_m^Ty_m
\end{pmatrix}	 \nonumber	\\
&=&
\begin{pmatrix}
\left (
\Xi_{m-1}^T \Pi_{m-1}
\Omega_{m-1} \Pi_{m-1}^T
\Xi_{m-1}
\right )
& \Xi_{m-1}^T\Pi_{m-1}\Psi_{m-1}^Ts_m \\
s_m^T\Psi_{m-1}\Pi_{m-1}^T\Xi_{m-1} & -s_m^TB_0s_m + s_m^Ty_m
\end{pmatrix} \nonumber	 \\	 
&=&
\begin{pmatrix}
\Xi_{m-1}^T\Pi_{m-1}^T & 0 \\
0 & 1 
\end{pmatrix}
\begin{pmatrix}
\Omega_{m-1}
& \Psi_{m-1}^Ts_m \\
s_m^T\Psi_{m-1} & -s_m^TB_0s_m + s_m^Ty_m
\end{pmatrix}
\begin{pmatrix}
\Pi_{m-1}^T\Xi_{m-1} & 0 \\
0 & 1 
\end{pmatrix}\nonumber\\ 
&=&
%	\begin{pmatrix}
%		\Xi_{m-1}^T\Pi_{m-1}^T & \ \ 0 & 0 \\
%		0 & -1 & 1
%	\end{pmatrix}
\Xi_m^T
\begin{pmatrix}
\Omega_{m-1}
& -\Psi_{m-1}^Ts_m & 0 \\
-s_m^T\Psi_{m-1} & -s_m^TB_0s_m + \gamma_m & \gamma_m \\
0 & \gamma_m &  s_m^Ty_m + \gamma_m
\end{pmatrix}
\Xi_m,  \label{eqn-middle}
\end{eqnarray} 
where
$$
	\Xi_m \defined
	\begin{pmatrix}
		\Pi_{m-1}^T\Xi_{m-1} & 0 \\
		0 & E_m
	\end{pmatrix},
	\quad \text{and} \quad
	E_m \defined
	\begin{pmatrix}
	-1 \\
	\ \ 1
	\end{pmatrix}.
$$
Note that the middle matrix in (\ref{eqn-middle}) can be expressed as
%\begin{equation}
%	\label{eq:PiOmegaPiT}
%	\begin{pmatrix}
%	\Omega_{m-1}
%	 & -\Psi_{m-1}^Ts_m & 0 \\
%	-s_m^T\Psi_{m-1}^T & -s_m^TB_0s_m + \gamma_m & \gamma_m \\
%	0 & \gamma_m &  s_m^Ty_m + \gamma_m
%	\end{pmatrix},
%\end{equation}
\begin{equation}
\label{eq:PiOmegaPiT}
\begin{pmatrix}
\Omega_{m-1}
& -\Psi_{m-1}^Ts_m & 0 \\
-s_m^T\Psi_{m-1} & -s_m^TB_0s_m + \gamma_m & \gamma_m \\
0 & \gamma_m &  s_m^Ty_m + \gamma_m
\end{pmatrix},
\end{equation} 
which is equivalent to
\begin{equation*}
	\begin{pmatrix}
		-S_{m-1}^TB_0S_{m-1} + \Gamma_{m-1}  & -L_{m-1}+\Gamma_{m-1} & -S_{m-1}^TB_0s_m & 0 \\
		-L_{m-1}^T+\Gamma_{m-1} & D_{m-1}+\Gamma_{m-1} & -Y_m^Ts_m  & 0 \\
		-s_m^TB_0S_{m-1} & -s_m^TY_m  & -s_m^TB_0s_m+\gamma_m & \gamma_m \\
		0  & 0 & \gamma_m   & y_m^Ts_m + \gamma_m
	\end{pmatrix}.
\end{equation*}
Substituting this into (\ref{eqn-middle}), yields
\begin{equation}\label{eqn-tildeHat}
\widetilde{\mathcal{M}}_m^{-1}=	
		\Xi_m^T
		\Pi_m
		\Omega_m
		\Pi_m^T
		\Xi_m,
\end{equation}
where $\Pi_m$ is defined in (\ref{eq:Pi_k}), replacing $k$ with
$m$, and $\Omega_m$ is defined in
(\ref{eqn-Omega}), replacing $m-1$ with $m$.  
Thus, 
$\widetilde{\mathcal{M}}_m = \widehat{M}_m$,
as defined in (\ref{eqn-G}).

We finish this case of the proof by showing that the last term in
\eqref{eqn-3x3} is equal to $\widehat{\Psi}_m \widehat{M}_m
\widehat{\Psi}_m^T$.  Substituting in \eqref{eq:GkPnew}
gives that the last term in (\ref{eqn-3x3})
can be written as
\begin{eqnarray*}
	\left ( \Psi_{m-1}\Pi_{m-1}^T \Xi_{m-1} \ \ B_0s_m \ \ y_m \right )
		\begin{pmatrix}
			I & \ \ 0 \\
			0 & -1 \\
			0 & \ \ 1
		\end{pmatrix}
	\widetilde{\mathcal{M}}_{m}
		\begin{pmatrix}
			I & \ \ 0 & 0 \\
			0 & -1 & 1
		\end{pmatrix}		
		\begin{pmatrix}
			\Xi_{m-1}^T \Pi_{m-1} \Psi_{m-1}^T \\
			(B_0s_m)^T \\
			y_m^T
		\end{pmatrix}.
\end{eqnarray*}
Using (\ref{eqn-tildeHat}) this simplifies to
$$	 \Psi_m \Pi_m^T
	 \begin{pmatrix}
	 	\Pi_{m-1}^T \Xi_{m-1} &  0 \\
		0 & E_m
	 \end{pmatrix}
	\left (
		\Xi_m^T
		\Pi_m
		\Omega_m
		\Pi_m^T
		\Xi_m
	\right )^{-1}
	 \begin{pmatrix}
	 	 \Xi_{m-1}^T\Pi_{m-1}& 0  \\
		0 & E_m
	 \end{pmatrix}
	 \Pi_m\Psi_m^T,$$
or, in other words,
%Simplifying this gives that the last term in (\ref{eqn-3x3}) is
$$
	 \Psi_m \Pi_m^T \Xi_m
	\left (
		\Xi_m^T
		\Pi_m
		\Omega_m
		\Pi_m^T
		\Xi_m
	\right )^{-1}
	\Xi_m^T \Pi_m \Psi_m^T,
$$
which is exactly $\widehat{\Psi}_m \widehat{M}_m
\widehat{\Psi}_m^T$.
Thus, for $\phi_m = \phi_m^{\SR1}$, the inductive step is proven.

\bigskip

\noindent \textbf{Case (ii):} We consider the case that $\phi_m \ne
\phi_m^{\SR1}$.  We begin by showing that $\widehat{M}_m=\mathcal{M}_m$,
given in \eqref{eq:Mbar}.  By Lemma 1, $\alpha_m + \beta_m \ne 0$.  Second,
$E_m = I_2$ (see \eqref{eq:Xi_k}), and $\gamma_m = \phi_m/(\alpha_m +
\beta_m)$ is well-defined (see \eqref{eq:Gamma}).  Then, the inverse of
$\mathcal{M}_m$ can be computed using arguments similar to those found
in~\cite{ErwM15}:
\begin{eqnarray}
		\mathcal{M}_m^{-1}
		&=&
		\begin{pmatrix}
		\widehat{M}_{m-1}^{-1} & -\widehat{M}_{m-1}^{-1}p_m & 0 \\
		-p_m^T\widehat{M}_{m-1}^{-1} & p_m^T\widehat{M}_{m-1}^{-1}p_m  + \tilde{\alpha}_m & \tilde{\beta}_m \\
		0 & \tilde{\beta}_m & \tilde{\delta}_m
		\end{pmatrix}, \label{eq:Ginv}
\end{eqnarray}
where 
\begin{equation}\label{eqn-bar}\tilde{\alpha}_m = \frac{\delta_m}{\alpha_m \delta_m - \beta_m^2}, 
\quad \tilde{\beta}_m=-\frac{\beta_m}{\alpha_m \delta_m -
\beta_m^2}\quad \text{and} \quad \tilde{\delta}_m = 
\frac{\alpha_m}{\alpha_m \delta_m - \beta_m^2}.\end{equation}
Simplifying the expressions in (\ref{eqn-bar}), yields
\begin{equation}\label{eqn-bar2}
\tilde{\alpha}_m = 
-s_m^TB_ms_m + \gamma_m, \quad \tilde{\beta}_m = \gamma_m, \quad
\tilde{\delta}_m= y_m^Ts_m + \gamma_m.
\end{equation}
We now simplify the entries of (\ref{eq:Ginv}) using the same approach
as in~\cite{ErwM15}.
Since $p_m = \widehat{M}_{m-1}\widehat{\Psi}_{m-1}^Ts_m$, then $\widehat{M}_{m-1}^{-1}p_m = \widehat{\Psi}_{m-1}^Ts_m$, 
giving us an expression for the (1,2) and (2,1) entries.  The (2,2) block entry 
is simplified by first multiplying (\ref{eqn-block22}) by $s_m^T$ on the left to obtain $s_m^TB_ms_m = s_m^TB_0s_m + p_m^T\widehat{M}_{m-1}^{-1}p_m$.  Then,
\begin{eqnarray*}
	p_m^T\widehat{M}_{m-1}^{-1}p_m + \tilde{\alpha}_m =  -s_m^TB_0s_m + s_m^TB_ms_m + \tilde{\alpha}_m 
	= -s_m^TB_0s_m + \gamma_m.
\end{eqnarray*}
Thus, using \eqref{eq:PiOmegaPiT}, (\ref{eq:Ginv}) can be written as
%\begin{eqnarray}
%	\mathcal{M}_m^{-1} 
%	&=& \nonumber
%		\begin{pmatrix}
%		\widehat{M}_{m-1}^{-1} & -\widehat{\Psi}_{m-1}^Ts_m & 0 \\
%		-s_m^T\widehat{\Psi}_{m-1} & -s_m^TB_0s_m\! +\! \gamma_m  
%		& \gamma_m \\
%		0 & \gamma_m
%		& y_m^Ts_m\! +\! \gamma_m
%		\end{pmatrix} \\
%	&=& \nonumber
%		\begin{pmatrix}
%		\Xi_{m-1}^T\Pi_{m-1}\Omega_{m-1}\Pi_{m-1}^T\Xi_{m-1}  & -\Xi_{m-1}^T\Pi_{m-1}\Psi_{m-1}^Ts_m & 0 \\
%		-s_m^T\Psi_{m-1}\Pi_{m-1}^T\Xi_{m-1} & -s_m^TB_0s_m\! +\! \gamma_m  
%		& \gamma_m \\
%		0 & \gamma_m
%		& y_m^Ts_m\! +\! \gamma_m
%		\end{pmatrix} \\
%	&=& \nonumber
%		\begin{pmatrix}
%		\Xi_{m-1}^T\Pi_{m-1}^T & 0 \\
%		0 & E_m
%		\end{pmatrix}
%		\begin{pmatrix}
%		\Omega_{m-1}
%		 & \Psi_{m-1}^Ts_m & 0  \\
%		s_m^T\Psi_{m-1}^T & -s_m^TB_0s_m + \gamma_m&
%		\gamma_m \\
%		 0 & \gamma_m &  y_m^Ts_m  + \gamma_m
%		\end{pmatrix}
%		\begin{pmatrix}
%		\Pi_{m-1}^T\Xi_{m-1} & 0 \\
%		0 & E_m 
%		\end{pmatrix}\\
%	&=& \nonumber
%		\Xi_m^T
%		\begin{pmatrix}
%		\Omega_{m-1}
%		 & \Psi_{m-1}^Ts_m & 0  \\
%		s_m^T\Psi_{m-1}^T & -s_m^TB_0s_m + \gamma_m&
%		\gamma_m \\
%		 0 & \gamma_m &  y_m^Ts_m  + \gamma_m
%		\end{pmatrix}
%		\Xi_m\\
%	&=&\label{eq:Mminv2}
%		\Xi_m^T \Pi_m \Omega_m \Pi_m^T \Xi_m,
%\end{eqnarray}
\begin{eqnarray}
\mathcal{M}_m^{-1} 
&=& \nonumber
\begin{pmatrix}
\widehat{M}_{m-1}^{-1} & -\widehat{\Psi}_{m-1}^Ts_m & 0 \\
-s_m^T\widehat{\Psi}_{m-1} & -s_m^TB_0s_m\! +\! \gamma_m  
& \gamma_m \\
0 & \gamma_m
& y_m^Ts_m\! +\! \gamma_m
\end{pmatrix} \\
&=& \nonumber
\begin{pmatrix}
\Xi_{m-1}^T\Pi_{m-1}\Omega_{m-1}\Pi_{m-1}^T\Xi_{m-1}  & -\Xi_{m-1}^T\Pi_{m-1}\Psi_{m-1}^Ts_m & 0 \\
-s_m^T\Psi_{m-1}\Pi_{m-1}^T\Xi_{m-1} & -s_m^TB_0s_m\! +\! \gamma_m  
& \gamma_m \\
0 & \gamma_m
& y_m^Ts_m\! +\! \gamma_m
\end{pmatrix} \\
&=& \nonumber
\begin{pmatrix}
\Xi_{m-1}^T\Pi_{m-1}^T & 0 \\
0 & E_m
\end{pmatrix}
\begin{pmatrix}
\Omega_{m-1}
& \Psi_{m-1}^Ts_m & 0  \\
s_m^T\Psi_{m-1} & -s_m^TB_0s_m + \gamma_m&
\gamma_m \\
0 & \gamma_m &  y_m^Ts_m  + \gamma_m
\end{pmatrix}
\begin{pmatrix}
\Pi_{m-1}^T\Xi_{m-1} & 0 \\
0 & E_m 
\end{pmatrix}\\
&=& \nonumber
\Xi_m^T
\begin{pmatrix}
\Omega_{m-1}
& \Psi_{m-1}^Ts_m & 0  \\
s_m^T\Psi_{m-1} & -s_m^TB_0s_m + \gamma_m&
\gamma_m \\
0 & \gamma_m &  y_m^Ts_m  + \gamma_m
\end{pmatrix}
\Xi_m\\
&=&\label{eq:Mminv2}
\Xi_m^T \Pi_m \Omega_m \Pi_m^T \Xi_m,
\end{eqnarray}
proving that $\widehat{M}_m=\mathcal{M}_m$.
Finally, using arguments similar to those in case (i), it can be shown that $$
	B_{m+1} = B_0 + \Psi_m \Pi_m^T \Xi_m (\Xi_m^T\Pi_m \Omega_m \Pi_m^T \Xi_m)^{-1} \Xi_m^T \Pi_m \Psi_m^T
	= B_0 + \widehat{\Psi}_m \widehat{M}_m \widehat{\Psi}_m^T,
$$
as desired.
 $\square$

\medskip

 There are two main differences in the compact representation for the full
 Broyden class (Theorem 1) and the restricted Broyden class (Section 2).
 First, in Theorem 1, $\Xi_k$ will always be the identity matrix for
 updates belonging to the restricted Broyden class.  Second, in
 (\ref{eq:compact}), the permutation matrices in the definitions of
 $\widehat{M}_k$ and $\widehat{\Psi}_k$, (equations (\ref{eqn-G}) and
 (\ref{eq:Phihat}), respectively) always cancel out in the restricted
 Broyden case. To emphasize that the permutation matrices do not cancel out
 for the general Broyden class updates, we use the notation $\widehat{M}_k$
 and $\widehat{\Psi}_k$, in lieu of $M_k$ and $\Psi_k$ as in the restricted
 Broyden case.

 Finally, we provide some insight regarding the permutation matrices
 (\ref{eq:Pi_k}).  The permutation matrix $\Pi_k$ acts in the following
 manner:
\begin{equation*}
	\Psi_k \Pi_k^T 
	=  \bigg ( B_0s_0 \ \ \cdots \ \ B_0s_{k-1} \ \ y_0 \ \ \cdots \ \ y_{k-1} \ \ B_0s_k \ \ y_k \bigg )
	= \bigg ( \Psi_{k-1} \ \ B_0s_k \ \ y_k \bigg ),
\end{equation*}
so that 
\begin{eqnarray}
	\Psi_k \Pi_k^T \Xi_k  &=&
	 \bigg ( \Psi_{k-1} \ \ B_0s_k \ \ y_k  \bigg )
	\begin{pmatrix}
		\Pi_{k-1}^T\Xi_{k-1} & 0 \\
		0 & E_{k}
	\end{pmatrix} 
	\nonumber \\
	&=&
	\bigg ( \Psi_{k-1}\Pi_{k-1}^T\Xi_{k-1} \ \ \  \big ( B_0s_k \ \ y_k \big ) E_k  \bigg )
	\nonumber \\
	&\vdots&
	\nonumber \\
	&=& \bigg ( \big (  B_0s_0 \ \ y_0 \big ) E_0 \ \ \  \big ( B_0s_1 \ \ y_1 \big ) E_1 \ \ \cdots \ \  
		\big ( B_0s_k \ \ y_k \big ) E_k \bigg ).
	\label{eq:PsiPiTXi}
\end{eqnarray}
In other words, when applied on the right of $\Psi_k$, the product
$\Pi_k^T\Xi_k$ permutes the columns of $\Psi_k$ and, using
the matrices $\{E_i\}$, combines columns of $\Psi_k$ whenever the
update is a rank-one update.

\bigskip

Unfortunately, computing $\widehat{M}_k$ is not straightforward.  In
particular, the diagonal matrix $\Gamma_k$ in Eq.\ \eqref{eq:Gamma}
involves $s_i^TB_is_i$ for each $i\in\{0, \ldots, k\}$, which requires
$B_i$ for $0 \le i \le k$.  In the next section, we propose a recursive method 
for computing $\widehat{M}_k$ that does not require storing the matrices
$B_i$ for $0 \le i \le k$.

\bigskip

\subsection{Computing $\widehat{M}_k$}
In this section, we 
propose a recursive method for computing $\widehat{M}_k$ from $\widehat{M}_{k-1}$.  
This method is based on the method proposed in~\cite{ErwM_CptInv}
for solving a linear system whose system matrix is
generated using the 
restricted Broyden class of updates.
In the proof of Theorem 1, we showed that
\begin{equation}\label{eq:widehatMk}
	\widehat{M}_k \ = \ 
	\begin{cases}
	\begin{pmatrix}
			\widehat{M}_{k-1} - \beta_kp_kp_k^T & \beta_k p_k \\
			\beta_k p_k^T & -\beta_k 
	\end{pmatrix}
	& \text{if $\phi_{k} = \phi_k^{\SR1}$} \\
	\begin{pmatrix}
			\widehat{M}_{k-1}+ \alpha_k p_kp_k^T & \alpha_k p_k & \beta_k p_k \\
			\alpha_k p_k^T & \alpha_k & \beta_k \\
			\beta_k  p_k^T & \beta_k & \delta_k
	\end{pmatrix}
	& \text{otherwise},
	\end{cases}
\end{equation}
which are given in  (\ref{eq:GkP}) and
(\ref{eqn-Mtildem-i}).  We now relate some of the entries in
$\widehat{M}_k$ with other stored or computable quantities involving
the pairs $\{s_i,y_i\}, i=0,\ldots,k$.
The vector $p_k$ can be computed as
\begin{equation}\label{eqn-formingp}
	p_k = \widehat{M}_{k-1}\widehat{\Psi}_{k-1}^Ts_k 
%	=
%	\widehat{M}_{k-1}\Xi_{k-1}^T\Pi_{k-1}
%	\begin{pmatrix}
%		(B_0S_{k-1})^T \\
%		Y_{k-1}^T
%	\end{pmatrix}
%	s_k
	=
	\widehat{M}_{k-1}
	\Xi_{k-1}^T\Pi_{k-1}
	\begin{pmatrix}
		S_{k-1}^TB_0s_k \\
		Y_{k-1}^Ts_k
	\end{pmatrix}.
\end{equation}
Note that in (\ref{eqn-formingp}), the vector $S_{k-1}^TB_0s_k$ is the first $k-1$ entries in the last
column of $S_k^TB_0S_k$, and the vector $Y_{k-1}^Ts_k$ is the first $k$
entries in the last column of
$Y_k^TS_k$.  Moreover, the entry $\alpha_k$,
given by $\alpha_k = -(1-\phi_k)/s_k^TB_ks_k$, can be
computed from the following:
\begin{equation}\label{eqn-makingalpha}
	s_k^TB_ks_k = s_k^T \bigg (B_0 + \widehat{\Psi}_{k-1}\widehat{M}_{k-1}\widehat{\Psi}_{k-1}^T \bigg ) s_k
	= s_k^TB_0s_k + s_k^T\widehat{\Psi}_{k-1}p_k.
\end{equation}
In (\ref{eqn-makingalpha}), the quantity $s_k^TB_0s_k$ is the $k$th 
diagonal entry in $S_k^TB_0S_k$, and $s_k^T\widehat{\Psi}_{k-1}p_k$ is the
inner product of $p_k$ and $\widehat{\Psi}_{k-1}^Ts_k$, the latter vector
already having been computed in (\ref{eqn-formingp}).  Recall that the
entry $\beta_k$ is given by $\beta_k = -\phi_k/y_k^Ts_k$, where $y_k^Ts_k$
is the $(k+1)$st diagonal entry in $S_k^TY_k$.  Finally, $\delta_k = ( 1 +
\phi_k s_k^TB_ks_k/y_k^Ts_k) / y_k^Ts_k $, which uses the previously
computed quantities $s_k^TB_ks_k$ and $y_k^Ts_k$.

For the initialization of $\widehat{M}_0$, notice 
that $\widehat{M}_0$ in (\ref{eqn-widehatM1}) can
be written as
\begin{equation}\label{eq:widehatM0}
	\widehat{M}_0 =
	\begin{cases}
	 -\beta_0
	& \text{if $\phi_0 = \phi_0^{\SR1}$}
	\\
	\begin{pmatrix}
		 \alpha_0 & \beta_0 \\
		\beta_0 & \delta_0
	\end{pmatrix}
	& \text{otherwise,}
	\end{cases}
\end{equation}
where $\alpha_0, \beta_0,$ and $\delta_0$ are defined as in (\ref{eqn-alphabetadelta}).

\medskip

In Algorithm \ref{alg:Bk+1}, we use
the recursions described above to compute $\widehat{M}_k$ given in (\ref{eq:widehatMk}).

\begin{algorithm}[!h]
\caption{This algorithm computes $\widehat{M}_k$ in (\ref{eq:widehatMk}).} \label{alg:Bk+1}
\SetAlgoNoLine
\KwIn{An initial $\phi_0$ and $B_0$;}
%\noindent \textbf{Algorithm 1.} \emph{This algorithm computes $\widehat{M}_k$ in (\ref{eq:widehatMk})}\\
Define $\widehat{M}_0$ using (\ref{eq:widehatM0}); \\
Define $\Psi_0 = (B_0s_0 \,\, y_0 )$;\\
\tab \FOR\  $j = 1:k$ 
\tabb $\widehat{\Psi}_{j-1}^Ts_j \gets \Xi_{j-1}^T \Pi_{j-1}\Psi_{j-1}^Ts_j$;
\tabb $p_j \gets \widehat{M}_{j-1} (\widehat{\Psi}_{j-1}^Ts_j$);
\tabb  $s_j^TB_js_j \gets s_j^TB_0s_j + ( s_j^T \widehat{\Psi}_{j-1})p_j$;
\tabb $\alpha_j \gets -(1-\phi_j)/ (s_j^TB_js_j$);
\tabb $\beta_j \gets -\phi_j/ (y_j^Ts_j$);
\tabb  $\delta_j \gets  ( 1 + \phi_j (s_j^TB_js_j)/(y_j^Ts_j) ) / (y_j^Ts_j)$;
\tabb Form $\widehat{M}_j$ using \eqref{eq:widehatMk};
\tab \END\ \\
\end{algorithm}

\noindent Note that the matrices $\Pi_{j-1}$ and $\Xi_{j-1}$ are not
explicitly formed in Algorithm \ref{alg:Bk+1}.  Instead, \eqref{eq:PsiPiTXi} can be used
to compute $\widehat{\Psi}_{j-1}$ in line 4 of Algorithm \ref{alg:Bk+1}.

\bigskip

\vspace{-6pt}
\section{Solving linear systems}
\vspace{-2pt}

Given the compact representation of $B_{k+1}$, we can solve 
\begin{equation}\label{eq:Br=z}
	B_{k+1} r = z,
\end{equation}
where $r, z \in \Re^n$, by computing the compact representation of the
inverse of $B_{k+1}$.  Intuitively speaking, computing the compact
representation of the inverse is due to the fact that $H_{k+1}\defined B_{k+1}^{-1}$ can also be written using a recursion
relation~\cite{DenM77}:
\begin{equation}\label{eq:H}
	H_{k+1}  = H_k  + \frac{1}{s_k^Ty_k}s_ks_k^T 
				- \frac{1}{y_k^TH_k y_k}H_k y_ky_k^TH_k 
				+ \Phi_k (y_k^TH_k y_k)v_kv_k^T,
\end{equation}
where $H_k\defined B_k^{-1}$,
$v_k = s_k/(y_k^Ts_k) - (H_k y_k)/(y_k^TH_k y_k)$,
and
\begin{equation}\label{eq:Phi}
	\Phi_k 
	= 
	\frac{(1-\phi_k)(y_k^Ts_k)^2}{(1-\phi_k)(y_k^Ts_k)^2 + \phi_k(y_k^TH_k y_k)(s_k^TB_k s_k)}.
\end{equation}
Note that when $\phi_k = \phi_k^{\SR1}$, then the corresponding $\Phi_k$ is given by
$$
	\Phi_k^{\SR1} = \frac{y_k^Ts_k}{y_k^Ts_k - y_k^TH_ky_k}.
$$

 In this section, we derive the compact
representation of the inverse of a Broyden class member.  This derivation
is similar to the process of finding the inverse of a member of the
restricted Broyden class presented in~\cite{ErwM_CptInv}.
 
Applying the Sherman-Morrison-Woodbury formula (see, e.g.,~\cite{GolV96})
to the compact representation of $B_{k+1}$ given in (\ref{eq:compact}),
gives that
$$
	B_{k+1}^{-1}= 
	B_0^{-1} + B_0^{-1} \widehat{\Psi}_k \left (-\widehat{M}_k^{-1} - \widehat{\Psi}_k^TB_0^{-1} \widehat{\Psi}_k \right )^{-1} \widehat{\Psi}_k^T B_0^{-1}.
$$
For quasi-Newton matrices it is conventional to let $H_i$ denote the inverse of $B_i$ for each $i$;
with this notation, the inverse of $B_{k+1}^{-1}$ is given by
\begin{equation}\label{eq:Hk+1}
	H_{k+1}= 
	H_0 + H_0 \widehat{\Psi}_k \left (-\widehat{M}_k^{-1} - \widehat{\Psi}_k^TH_0 \widehat{\Psi}_k \right )^{-1} \widehat{\Psi}_k^T H_0.
\end{equation}
Using the definition of  $\widehat{\Psi}_k$ in (\ref{eq:Phihat}) gives that
\begin{equation*}
	\widehat{\Psi}_k^TH_0 \widehat{\Psi}_k =
	\Xi_k^T\Pi_k
	\begin{pmatrix}
		S_k^TB_0S_k & S_k^TY_k \\
		Y_k^TS_k & Y_k^TH_0Y_k
	\end{pmatrix}
	\Pi_k^T\Xi_k,
\end{equation*}
and thus,
$$
	-\widehat{M}_k^{-1} - \widehat{\Psi}_k^TH_0 \widehat{\Psi}_k = 
	-\Xi_k^T\Pi_k
	\begin{pmatrix}
		\Gamma_k & D_k+R_k+\Gamma_k \\
		D_k + R_k^T + \Gamma_k & D_k + \Gamma_k + Y_k^TH_0Y_k
	\end{pmatrix}
	\Pi_k^T\Xi_k.
$$
Substituting $H_0 \widehat{\Psi}_k = H_0 (B_0S_k  \ Y_k )\Pi_k^T \Xi_k = (S_k  \ H_0Y_k )\Pi_k^T\Xi_k$ into  (\ref{eq:Hk+1}) gives 
the compact representation for the inverse of any member of the full Broyden class:
\begin{equation}\label{eq:Hk+1compact}
	H_{k+1} = H_0 + \wPsi_k \wM_k \wPsi_k^T,
\end{equation}
where $\wPsi_k = (S_k  \ H_0Y_k )\Pi_k^T\Xi_k$ and
\begin{equation}\label{eq:Mtilde}
	\wM_k \equiv
	\left (
	-\Xi_k^T\Pi_k
	\begin{pmatrix}
		\Gamma_k & D_k+R_k+\Gamma_k \\
		D_k + R_k^T + \Gamma_k & D_k + \Gamma_k + Y_k^TH_0Y_k
	\end{pmatrix}
	\Pi_k^T\Xi_k
	\right )^{-1}.
\end{equation}
\\
\noindent \textbf{Computing $\wM_k$.}
Using an approach similar to how $\widehat{M}_k$ is computed, 
$\wM_k$ can be computed as follows:
\begin{equation}\label{eq:widetildeM}
\wM_k = 
	\begin{cases}
	\begin{pmatrix}
		\wM_{k-1} - \tilde{\beta}_k\tilde{p}_k\tilde{p}_k^T & -\tilde{\beta}_k\tilde{p}_k \\
		-\tilde{\beta}_k\tilde{p}_k^T & -\tilde{\beta}_k
	\end{pmatrix}
	& \text{if $\Phi_k = \Phi_k^{\SR1}$}
	\\
	\begin{pmatrix}
		\wM_{k-1}
			+\tilde{\delta}_k\tilde{p}_k\tilde{p}_k^T & \tilde{\beta}_k\tilde{p}_k 
			& \tilde{\delta}_k \tilde{p}_k \\
		\tilde{\beta}_k\tilde{p}_k^T & \tilde{\alpha}_k & \tilde{\beta}_k \\
		\tilde{\delta}_k\tilde{p}_k^T &\tilde{\beta}_k & \tilde{\delta}_k
	\end{pmatrix}
	& \text{otherwise},
	\end{cases}
\end{equation}
where
\begin{equation}\label{eq:tildealphabetadelta}
	\tilde{\alpha}_k = \frac{1}{s_k^Ty_k}+ \Phi_k\frac{y_k^TH_ky_k}{(s_k^Ty_k)^2}, \quad
	\tilde{\beta}_k = -\frac{\Phi_k}{y_k^Ts_k}, \quad
	\tilde{\delta}_k = -\frac{1-\Phi_k}{y_k^TH_ky_k},
\end{equation}
and $\tilde{p}_k = \wM_{k-1}\wPsi_{k-1}^Ty_k$. The initial matrix $\wM_0$ is given by the following:
\begin{equation}\label{eq:widetildeM0}
	\wM_0 =
	\begin{cases}
	 -\tilde{\beta}_0
	& \text{if $\Phi_0 = \Phi_0^{\SR1}$}
	\\
	\begin{pmatrix}
		 \tilde{\alpha}_0 & \tilde{\beta}_0 \\
		\tilde{\beta}_0 & \tilde{\delta}_0
	\end{pmatrix}
	& \text{otherwise},
	\end{cases}
\end{equation}
where $\tilde{\alpha}_0, \tilde{\beta}_0,$ and $\tilde{\delta}_0$ are defined as in 
\eqref{eq:tildealphabetadelta} with $k = 0$.  A practical iterative method
to solve equations of the form (\ref{eq:Br=z}) is given in Algorithm \ref{alg:Br=z}.

\begin{algorithm}[!h]
\caption{This algorithm solves $B_{k+1}r = z$.}\label{alg:Br=z}
\SetAlgoNoLine
\KwIn{An initial $\phi_0$, $B_0$, and $H_0$;}
Define $\widehat{M}_0$ using \eqref{eq:G0} and $\wM_0 = $ using \eqref{eq:widetildeM0};
\tab \FOR\ $j = 1:k$
\tabb Compute $s_j^TB_js_j$ using Algorithm \ref{alg:Bk+1};
\tabb $\wPsi_{j-1}^Ty_j \gets  \Xi_{j-1}^T \Pi_{j-1} \Psi_{j-1}^TH_0 y_j$;
\tabb $\tilde{p}_j \gets \wM_{j-1} (\wPsi_{j-1}^Ty_j$);
\tabb $y_j^TH_jy_j \gets y_j^TH_0y_j + ( y_j^T \wPsi_{j-1})\tilde{p}_j$;
\tabb $\Phi_j \gets (1 - \phi_j)(y_j^Ts_j)^2 / ((1 - \phi_j)(y_j^Ts_j)^2 + \phi_j(y_j^TH_jy_j)(s_j^TB_js_j))$;
\tabb $\tilde{\alpha}_j \gets (1+\Phi_j(y_j^TH_jy_j)/(y_j^Ts_j)) / (y_j^Ts_j)$;
\tabb $\tilde{\beta}_j \gets -\Phi_j/ (y_j^Ts_j$);
\tabb $\tilde{\delta}_j \gets -(1 - \Phi_j) / (y_j^TH_jy_j)$;
\tabb Form $\wM_j$ using \eqref{eq:widetildeM};
\tab \END\ \\
$\wPsi_{k} \gets  H_0 \Psi_{k} \Pi_{k}^T\Xi_{k}$;\\
$r = H_0z + \wPsi_k\wM_k \wPsi_k^Tz$;\\
\end{algorithm}

\bigskip

%%%%%%%%%%%%%%%%%%%%%%%%%%%%%%%%%%%%%%%%%%%%%%%%%%%%%%%%%%%%%%%%%%%
\vspace{-6pt}
\section{Numerical experiments}
\vspace{-2pt}

In this section we test the accuracy of Algorithm \ref{alg:Bk+1} to compute the
  compact representation by comparing it with the matrix obtained using
the Broyden update formula \eqref{eqn-1param}.  In addition, we demonstrate
that solves with $B_{k+1}$ in \eqref{eq:Br=z} can be done efficiently using
Algorithm \ref{alg:Br=z} with respect to both accuracy and time.  For these
experiments, we used five
(limited-memory) quasi-Newton pairs to compute $B_{k+1}$.  To
generate quasi-Newton pairs, we simulated a line-search method where the
iterates are updated as follows:
$$
	x_{j+1} = x_j - \alpha_j B_j^{-1}g_j, \quad \text{for $1 \le j \le 4$},
$$	
where $\alpha_j \in [0, 1]$ was generated randomly.  To initialize the
process, we randomly generated initial points $x_0$ and $x_1$ so that
  $s_0 = x_1 - x_0$.  The corresponding gradients, $g_j = \nabla f(x_j)$
  for $0 \le j \le 5$, were also generated randomly in order to form
 $y_j = g_{j+1} - g_j$ for $0 \le j \le 4$.
  The matrix $B_0$ was initially defined as
$B_0 = \gamma I$, where $\gamma > 0$ was randomly generated.
We considered four experiments where we vary the value of $\phi_i$ at each
iteration $i$.  In particular, we  chose values of $\phi_i$ according
to the scheme given in Table \ref{table:phi}.
We ran each experiment ten times with $n = 10, 100,$ and $1,000$ and report results.
\begin{table}[!h]
\centering
\begin{tabular}{|c|c|c|c|c|c|}%C{1.3cm}|C{1.3cm}|C{1.8cm}|C{1.3cm}|C{1.3cm}|}
\hline
Experiment &$\phi_0$   &  $\phi_1$  \quad  & $\phi_2$ & $\phi_3$ & $\phi_4$\\ \hline
1 & $\phi_0 < 0$ & 1 & $0 < \phi_2 < 1$ & 0 & $\phi_4 > 1$ \\ \hline
2 & $\phi_0 < 0$ & 1 & $\phi_2^{\SR1}$ & 0 & $\phi_4 > 1$ \\ \hline
3 & $\phi_0 < 0$ & 1 & $\phi_2^{\SR1}$ & $\phi_3^{\SR1}$ & $\phi_4 > 1$ \\ \hline 
4 & $\phi_0^{\SR1}$ & 1 & $\phi_2^{\SR1}$ & 0 & $\phi_4 > 1$ \\ \hline
\end{tabular}
\caption{The values of $\phi_i$ for $0 \le i \le 4$ for each experiment.  The choice of $\phi_1 = 1$ corresponds to 
the BFGS update while $\phi_3 = 0$ corresponds to the DFP update.  Note that Experiment 1 does not use SR1 updates.}
\label{table:phi}
\end{table}

\subsection{Accuracy of the compact representation} 
To test the accuracy of the compact representation, we form each
  $B_{k+1}$ using (\ref{eq:compact}) together with the proposed compact
  formulation given in Theorem 1.  (In particular, we use Algorithm \ref{alg:Bk+1} to
  form $\widehat{M}_k$.)  We denote the resulting matrix by $B_{k+1}^{\text{CR}}$.
 In Table 2, we report the average relative
  error of the compact representation in the Frobenius norm:
$$
	\text{Relative error } = \frac{ \| B_{k+1} -
  	B_{k+1}^{\text{CR}}\|_F}{\|B_{k+1} \|_F}, 
$$ 
where $B_{k+1}$  is computed using (\ref{eqn-1param}).

\begin{table}[h]
\centering 
	\begin{tabular}{|c|c|c|c|c|c|}
		\hline
		$n$ & Exp.\ 1 & Exp.\ 2 & Exp.\ 3 & Exp.\ 4\\ \hline
		100 &  \texttt{1.1315e-13}& \texttt{1.3383e-11}& \texttt{1.6749e-12}& \texttt{2.2855e-14} \\ 
		1,000 & \texttt{3.2039e-14}& \texttt{1.1225e-14}& \texttt{5.4247e-15}& \texttt{1.0155e-15} \\ 
		10,000 &\texttt{1.3426e-13}& \texttt{8.5453e-14}& \texttt{1.9969e-13}& \texttt{2.8354e-16}\\ \hline
	\end{tabular}\\
		\caption{Average relative error over ten different trials for each experiment with $n = 100, 1,000,$ and $10,000$.}
		\label{table:MatrixRelErr}
\end{table}
The small relative errors in Table \ref{table:MatrixRelErr} reflects the fact that the
proposed compact representation for the full Broyden class of
quasi-Newton matrices is correct; moreover, the relative errors suggest
that Algorithm \ref{alg:Bk+1} provides a method to compute the compact representation
to high accuracy.

\subsection{Accuracy of  the compact representation of the inverse}
In these experiments, we test the accuracy of Algorithm \ref{alg:Br=z} to solve
linear systems of the form $B_{k+1} r =z$, where $r,z\in\Re^n$ and
$B_{k+1}$ is a quasi-Newton matrix.  The matrix $B_{k+1}$ is generated
using five quasi-Newton pairs as described in the beginning of this
section.  Moreover, the righthand side $z$ is randomly generated for each
experiment.  In Table \ref{table:RelResid}, we present the average residual error using the
two-norm:
$$
\text{Relative error} = \frac{ \|B_{k+1}r^{\text{ICR}}-z\|_2}{\|z\|_2},
$$ 
where $r^{\text{ICR}}$ is the solution to $B_{k+1}r=z$ using the inverse
compact representation computed by Algorithm \ref{alg:Br=z}.  These results suggest that
the compact representation of the inverse can be used to solve linear
systems to high accuracy.

\begin{table}[h]
	\centering
	\begin{tabular}{|c|c|c|c|c|c|}
		\hline
		$n$ & Exp.\ 1 & Exp.\ 2 & Exp.\ 3 & Exp.\ 4\\ \hline
		100 &  \texttt{4.0158e-13}& \texttt{1.342e-10}& \texttt{1.3065e-09}& \texttt{2.8160e-14} \\ 
		1,000 & \texttt{1.518e-14}& \texttt{7.6460e-14}& \texttt{6.1744e-14}& \texttt{1.8431e-13} \\ 
		10,000 &\texttt{2.4175e-12}& \texttt{1.6079e-12}& \texttt{4.3284e-12}& \texttt{1.8795e-14}\\ \hline
	\end{tabular}\\
	\caption{Average relative error over ten different trials for each experiment with $n = 100, 1,\!000,$ and $10,\!000$.}
	\label{table:RelResid}
\end{table}

In addition, during the experiments, the computational time of the
  proposed method was recorded and compared to a similar solve using the
  \MATLAB{} ``backslash''.  In particular, with the same quasi-Newton
  pairs, the backslash command was used to solve $B_{k+1}r=z$, where
  $B_{k+1}$ was formed using \eqref{eqn-1param}.  The times required were
  averaged for each experiment and for each value of $n$. These results are
  given in Table 4 and do not include the time \MATLAB{} required to
    form $B_{k+1}$. Note that the average computational times in
  Table \ref{table:time} indicate that as $n$ increases using Algorithm \ref{alg:Br=z} becomes
  significantly less computationally expensive than using the backslash
  command.

%
%\begin{table}[h]
%	\centering
%	\begin{tabular}{|c|c|c|c|c|c|c|}
%		\hline
%		$n$& & Exp.\ 1 & Exp.\ 2 & Exp.\ 3 & Exp.\ 4\\ \hline
%		100 & Time Compact&  \texttt{9.9143e-04}& \texttt{7.1645e-04}& \texttt{6.6531e-04}& \texttt{5.8214e-04} \\ 
%		&Time Backslash &\texttt{3.3208e-04} &\texttt{3.7331e-04} & \texttt{3.1091e-04}&\texttt{3.4517e-04} \\\hline
%		1,000 &Time Compact & \texttt{1.0707e-03}& \texttt{1.0768e-03}& \texttt{1.0551e-03}& \texttt{9.8053e-04} \\ 
%		&Time Backslash &\texttt{2.8957e-02} &\texttt{2.9188e-02} &\texttt{2.9644e-02} &\texttt{3.0362e-02} \\\hline
%		10,000 &Time Compact &\texttt{3.6523e-03}& \texttt{3.5308e-03}& \texttt{3.2106e-03}& \texttt{2.7792e-03}\\ 
%		& Time Backslash&\texttt{1.1335e+01} &\texttt{1.1368e+01} &\texttt{1.1341e+01} &\texttt{1.1429e+01} \\\hline
%	\end{tabular}\\
%	\caption{Average computational times for $n = 100, 1,000$, and $10,000$.}
%	\end{table}  

\begin{table}[h]
\centering
\begin{tabular}{|c|c|c|c|c|c|c|}
\cline{2-7}
\multicolumn{1}{c|}{} 
& \multicolumn{2}{c|}{$n = 100$} 
& \multicolumn{2}{c|}{$n = 1,000$} 
& \multicolumn{2}{c|}{$n = 10,000$} \\
\hline
Exp. & ICR & MATLAB
& ICR & MATLAB
& ICR & MATLAB
\\ \hline
1 & \texttt{9.9e-04} 
& \texttt{3.3e-04} 
& \texttt{1.1e-03}
& \texttt{2.9e-02}
& \texttt{3.7e-03}
& \texttt{1.1e+01} 
\\
2 & \texttt{7.2e-04} 
& \texttt{3.7e-04} 
& \texttt{1.1e-03}
& \texttt{2.9e-02}
& \texttt{3.5e-03}
& \texttt{1.1e+01}
\\
3 & \texttt{6.7e-04} 
& \texttt{3.1e-04} 
& \texttt{1.1e-03}
& \texttt{3.0e-02}
& \texttt{3.2e-03}
& \texttt{1.1e+01}
\\
4 & \texttt{5.8e-04} 
& \texttt{3.5e-04}
& \texttt{9.8e-04}
& \texttt{3.0e-02}
& \texttt{2.8e-03}
& \texttt{1.1e+01} 
\\
\hline
\end{tabular}
\caption{Average computational times for solving $B_{k+1}r = z$ using the inverse compact representation (ICR)
of the inverse and the MATLAB{} ``backslash'' command
with $n = 100, 1,000$, and $10,000$.}
\label{table:time}
\end{table}	
	
\vspace{-6pt}
\section{Conclusion}
\vspace{-2pt}

We derived the compact formulation for members of the full Broyden class of
quasi-Newton updates.  The compact representation allows for different
$\phi_k$ at each iteration as well as different ranks of updates.  With
this compact formulation, we demonstrated how to solve linear systems
defined by these limited-memory quasi-Newton matrices.  Numerical results
suggest that the compact representation can be computed to high
  accuracy and that we can solve \eqref{eq:Br=z} efficiently and
accurately using the compact representation of the inverse of $B_{k+1}$.  Future work includes
integrating this linear solver inside large-scale optimization methods.

\vspace{-6pt}
\section{Acknowledgments}
\vspace{-2pt}

The authors would like to thank Lasith Adhikari and Johannes Brust for helpful discussions regarding this work.
This research is supported by NSF grants CMMI-1334042 and CMMI-1333326.

\bibliographystyle{siam}
\bibliography{DeGuchyErwayMarcia}

\end{document}